\magnification=\magstep1
\input amstex
\documentstyle{amsppt}
\NoBlackBoxes \vsize=7in \hsize=5.5in

\def\mycal{\mathfont@\rsfs}
\csname rsfs \endcsname \topmatter
 
\def\sdp{\medspace \times \kern -1.9pt  
         \vrule width0.4pt height 4.7pt depth -0.3pt \medspace}
\hsize=5.2in
\vsize=6.8in
\centerline {\bf A RELATIVE VERSION OF CONNES`  $\chi(M)$ INVARIANT}
\centerline {\bf AND EXISTENCE OF ORBIT INEQUIVALENT ACTIONS} 
\vskip .1in
\centerline {\rm by}
\vskip .1in
\centerline {\rm ADRIAN IOANA}
\address Math Dept  
UCLA, Los Angeles, CA 90095-155505\endaddress
\email adiioana\@math.ucla.edu \endemail
\topmatter
\abstract We consider 
a new orbit equivalence invariant for measure-preserving actions 
of groups on the probability space, $\sigma:G\rightarrow$ Aut$(X,\mu)$
 , denoted  $\chi_0(\sigma;G)$ and 
defined as the "intersection" of the 1-cohomology group,
H$^1(\sigma,G)$, with Connes' invariant, $\chi(M)$,
of the cross product von Neumann algebra, $M=L^\infty(X,\mu)\rtimes_\sigma G$.
We calculate $\chi_0(\sigma;G)$ for certain actions of groups of the form $G=H\times K$ with $H$ 
non-amenable and $K$ infinite amenable and we deduce that any such group has
uncountably many orbit inequivalent actions.
 
\endabstract
\endtopmatter
\document
\vskip .1in
\head 0. Introduction. \endhead
\vskip .1in
\noindent
 Recall that two  measure preserving (m.p.) actions $\sigma_1,\sigma_2$  of two countable, discrete groups $G_1,G_2$,
on  standard probability measure spaces $(X_1,\mu_1),(X_2,\mu_2)$, 
are said to be {\it orbit equivalent} (OE)
if there is an isomorphism $f:(X_1,\mu_1)\rightarrow (X_2,\mu_2)$ such that $f(G_1x)=G_2(f(x))$ for $\mu_1$ 
a.e. $x\in X_1$, or, equivalently, if the induced equivalence relations, $\Cal R_{\sigma_1,G_1},\Cal R_{\sigma_2,G_2}$ are isomorphic.
 
Strikingly, it was proven that any  action $\sigma$ of any  amenable group $G$ induces a hyperfinite equivalence relation (unique, up to isomorphism), thus implying that   any two actions of any two amenable groups are OE. This result was first obtained by Dye([Dy]) in the case $G=\Bbb Z$  and  by Ornstein-Weiss([OW],[CFW]) in general.
The converse is also true: any non-amenable group has at least 2 non-OE actions([CW],[Sc1,2],[Hj]).
Thus, it is  natural to consider the problem of finding  groups which admit many non OE actions.
In this direction,  rigidity
phenomena were used to exhibit classes of groups possesing uncountably many
 non-OE actions:
property(T) groups([Hj]), products of hyperbolic groups([MoSh]), free groups([GaPo]) and weakly rigid 
groups([Po]).  

In this paper,   we introduce a relative version of Connes` $\chi(M)$ invariant ([Co1]) for Cartan subalgebra
inclusions $A\subset M$, denoted $\chi_0(M;A)$, which captures 
the "approximately $A$-inner, centrally trivial"  outer automorphisms of $M$. The $\chi_0$ invariant of  a free, ergodic, m.p.  action $\sigma:G\rightarrow\text{Aut}(X,\mu)$ is then defined as the $\chi_0$ invariant of the corresponding Cartan subalgebra inclusion, i.e. $\chi_0(\sigma;G)=\chi_0(L^{\infty}(X,\mu)\rtimes_{\sigma}G;L^{\infty}(X,\mu))$. Note that  $\chi_0(\sigma;G)$ is is an {\it orbit equivalence} invariant([FM]) and that it is  a subgroup of the 1-cohomology group H$^1(\sigma;G)$. 

\vskip 0.05in 
\proclaim{Theorem} Let $G=H\times K$, where $H$ is a non-amenable group and $K$ is an $\infty$ amenable group and let $\Gamma=\bigoplus_{i\geq 0}\Delta_i$, where $\Delta_i$ are finite groups. For a standard probability space $(X_0,\mu_0)$, let $\sigma$ be  
the generalized Bernoulli shift action of $G$ on $(X,\mu)=\prod_{g\in H\sqcup K}(X_0,\mu_0)$, given by the natural action of $G$ on the set $H\sqcup K$. Let $\beta$ be a free, m.p. action of $\Gamma$ on $(X_0,\mu_0)$ and define $\tilde\beta_{\gamma}=\otimes_{g\in H\sqcup K}(\beta_{\gamma})_g,\forall \gamma\in\Gamma$ to be the induced action on $(X,\mu)$.
 
If $\sigma^{\Gamma}$ denotes the restriction of $\sigma$ to the fixed point algebra $\{a\in L^{\infty}(X,\mu)|\tilde{\beta}_{\gamma}(a)=a,\forall\gamma\in\Gamma\}$, then $\sigma^{\Gamma}$ is a free, ergodic, m.p. action such that $\chi_0(\sigma^{\Gamma};G)=\text{Char}(\Gamma)$.
\endproclaim
\vskip 0.03in
Since for every set $\Cal P$ of prime numbers, we have that $\prod_{p\in\Cal P}\Bbb Z_p=\text{Char}(\bigoplus_{p\in \Cal P}\Bbb Z_p)$, we deduce the following:
\vskip 0.03in

\proclaim{Corollary} Any group of the form $G=H\times K$, where $H$ is a non-amenable group and $K$ is an $\infty$ amenable group, admits uncountably many non orbit equivalent actions.
\endproclaim
\vskip 0.05in

In Section 1 we give the definition of $\chi_0(\sigma,G)$  for an action $(\sigma,G)$ and we notice that it is a {\it stable orbit equivalence} invariant. Starting with Section 2, we  concentrate on the actions $(\sigma^{\Gamma},G=H\times K)$ as in the above Theorem and we give a first estimate of $\chi_0$. This reduces the calculation of $\chi_0$ to a concrete problem in terms of automorphisms of the hyperfinite II$_1$ factor, $R$.
In Section 3 we deal with this problem, showing that  for certain actions of groups $\Gamma$ as above  on the hyperfinite II$_1$ factor $R$, the following holds true: $\theta\in$ Aut($R$) acts trivially on the  sequences from $R^{\Gamma}$ which are central in $R$  iff $\theta\in\overline{\Gamma}$ (modulo inner automorphisms).
The final section is concerned with the proof of the above Theorem.

 \vskip 0.05in  
{\it Acknowledgement.} I would like to thank Professor Sorin Popa for proposing the original idea from which this work developed. I greatly benefitted from many discussions with him and I am grateful for his constant guidance and encouragement.
   
\vskip .1in
\head 1.Definition of $\chi_0.$ \endhead
\vskip .1in
\noindent 
  
In [Co1], A. Connes defined an invariant $\chi(M)$ for a II$_1$ factor $M$ and used computations of it to provide a II$_1$ factor non-antisomorphic to itself. 
 The purpose of this section is to define an  analogous invariant for inclusions $(A\subset M)$ of Cartan subalgebras into II$_1$ factors. Restricting to such inclusions that arise from measure preserving actions of countable, discrete groups on standard probability spaces we obtain an OE invariant, which we subsequently use to produce non-OE actions.

{\it 1.1. Groups of automorphisms.}
We begin by reviewing several definitions and notations. Given a  {\rm II}$\sb 1$ factor $M$ we denote by Aut$(M)$ its 
group of automorphisms endowed with the topology given by the pointwise norm  $\|.\|_2$ convergence(see [Co1,Co2]) and by Int$(M)$
 the subgroup of inner automorphisms of $M$. A bounded sequence $(x_n)_n\subset M$  is called {\it central} if for any element $y\in M$
 we have that $\lim_{n\rightarrow \infty}\|[x_n,y]\|_2=0$. Following Connes, an automorphism $\theta$ of $M$ is called 
{\it centrally trivial} if for any central sequence $(x_n)_n\subset M$  we have that
 $\lim_{n\rightarrow \infty}\|\theta(x_n)-x_n\|_2=0$. The group of centrally trivial automorphisms of $M$ is denoted Ct$(M)$([Co1,Co2]). Similarly, for an inclusion of von Neumann algebras $N\subset M$, we denote by Ct$(M,N)$ the group of automorphisms of $M$ that act trivially on the bounded sequences of $N$ which are central in $M$([Ka]).

 A powerful tool in the study of central sequences is the {\it ultrapower} algebra $M^{\omega}$ defined as the quotient  $l^{\infty}(\Bbb N,M)/{\Cal I_\omega}$, where $\Cal I_{\omega}=\{x=(x_n)_n\in l^{\infty}(\Bbb N;M)|\lim_{n\rightarrow\omega}||x_n||_2=0\}$ and $\omega$ is a {\it free ultrafilter} on $\Bbb N$. Any $\theta\in\text{Aut}(M)$ induces an automorphism $\theta^{\omega}$ of $M^{\omega}$ given by $\theta^{\omega}((x_n))=(\theta(x_n))_n$ and we have  that $\text{Ct}(M)=\{\theta\in\text{Aut}|\theta^{\omega}(x)=x,\forall x\in M'\cap M^{\omega}\}$ and similarly that $\text{Ct}(M,N)=\{\theta\in\text{Aut}(M)|\theta^{\omega}(x)=x,\forall x\in M'\cap N^{\omega}\}$.

\vskip 0.05in 
{\it 1.2. Automorphisms coming from 1-cocycles.}
For a  {\rm II}$\sb 1$ factor $M$ with a Cartan subalgebra $A$ we denote by Aut$\sb 0(M;A)$  the group of automorphisms of $M$ which leave $A$ pointwise fixed and by Int$\sb 0(M;A)$ 
 the subgroup of inner automorphisms implemented by the unitaries of $A$. Also, we consider the quotient group Out$\sb 0(M;A)$:= Aut$\sb 0(M;A)$/ Int$\sb 0(M;A)$ together with the quotient map $\varepsilon_0:$ 
Aut$_0(M;A)\rightarrow$ Out$_0(M;A)$.
 
Let  now $\sigma:G\rightarrow$ Aut$(X,\mu)$ be a free, ergodic, m.p. action on a standard probability space. Then $\sigma$ defines an integral preserving action of $G$ on the diffuse abelian von Neumann algebra $L^{\infty}(X,\mu)$; conversely, every integral preserving action on $L^{\infty}(X,\mu)$ gives rise to a m.p. action on $(X,\mu)$. Therefore, we can unambiguously identify the two actions.
 
Further, consider the associated group measure space factor $M=L^{\infty}(X,\mu)\rtimes_{\sigma}G$ together with its Cartan subalgebra $A=L^{\infty}(X,\mu)\subset M$([MvN]). Then  we have a canonical isomorphism
from  Aut$\sb 0(M;A)$ onto the group of 1-cocycles, Z$\sp 1(\sigma,G)=\{w:G\rightarrow\Cal U(L^{\infty}(X,\mu))|w_{gh}=w_g\sigma_g(w_h),\forall g,h\in G\}$, which carries Int$\sb 0(M;A)$ onto the subgroup of 1-coboundaries,
 B$\sp 1(\sigma,G)$(see [Si],[Po]). From this it follows that Out$_0(M;A)$ is isomorphic to the 1-cohomology group of $\sigma$, 
H$^1(\sigma,G)$:= Z$^1(\sigma,G)$/ B$^1(\sigma,G)$. We will also use the notations Int$_0(\sigma;G)$ and Aut$_0(\sigma;G)$ for the groups of automorphisms Int$_0(A\rtimes_{\sigma}G;A)$ and Aut$_0(A\rtimes_{\sigma}G;A)$, respectively.
\vskip 0.07in

Recall that by [FM], two free, m.p. actions $\sigma_i:G_i\rightarrow$ Aut$(A_i)$ are  {\it orbit equivalent} iff 
$(A_1\subset A_1\rtimes_{\sigma_1} G_1)\simeq(A_2\subset  A_2 \rtimes_{\sigma_2}
 G_2)$, i.e. if there exists an isomorphism $\theta:A_1\rtimes_{\sigma_1}G_1\rightarrow A_2\rtimes_{\sigma_2}G_2$ such that $\theta(A_1)=A_2$. Thus, any invariant of Cartan subalgebra inclusions naturally gives rise to an OE invariant for actions. 
 In particular, this entails that H$^1(\sigma,G)$ is an OE invariant.

For strongly ergodic actions $\sigma$, i.e. actions that have only trivial asymptotically invariant sequences, H$^1(\sigma;G)$ is a Polish group ([Sc1]),
which is  countable, discrete whenever $G$ is weakly rigid([Sc1],[Po]).
As shown by S.Popa([Po]), for w-rigid groups $G$ and for actions $\sigma$ having 
good deformation properties, the 1-cohomology group, H$^1(\sigma;G)$ is calculable and can be any countable abelian group.
Thus, in the case of w-rigid groups, computations of H$^1(\sigma;G)$ can be used to provide uncountably many 
non-OE actions.
 
At the other 
end, if the action is not strongly ergodic (equivalently, if Int$\sb 0(M;A)$ is not closed in the topology inherited from
Aut$(M)$) then this invariant 
is a large, non-Polish group which can thus not be used to distinguish between orbit inequivalent
actions.
Instead, the following  subgroup of Out$_0(M;A)\simeq$ H$^1(\sigma,G)$ can be used:
\vskip 0.05in 
{\it 1.3. Definition.} Let $M$ be a separable {\rm II}$\sb 1$ factor with a Cartan subalgebra $A$.
We define $\chi_0(M;A)= {\overline {\text{Int}}}_0(M;A)\cap \text{Ct}M/ {{\text{Int}}
_0(M;A)}$. 

\vskip 0.05in
In particular, if $\sigma:G\rightarrow$ Aut$(X,\mu)$ is a free ergodic m.p. action on a standard probability 
space then
 $\chi_0(\sigma,G):= \chi_0(L^{\infty}(X,\mu) \rtimes_\sigma G;L^{\infty}(X,\mu))$ defines an orbit equivalence invariant.
This invariant is a subgroup of H$^1(\sigma,G) \simeq$ Out$\sb 0(L^{\infty}(X,\mu)\rtimes_\sigma G;L^{\infty}(X,\mu))$, 
 which is trivial if $G$ is amenable([OW])  and also, by definition, if $\sigma$ is strongly ergodic.

Now, recall the following result from [FM]:  two free, ergodic, m.p. actions $\sigma_i:G_i\rightarrow$ Aut$(A_i)$ are {\it stably orbit equivalent} iff 
there exist two projections $ p_i\in A_i$, $i=1,2$ such that 
$(A_1p_1\subset p_1(A_1\rtimes_{\sigma_1} G_1)p_1)\simeq(A_2p_2\subset p_2( A_2 \rtimes_{\sigma_2}
 G_2)p_2)$. 
\vskip .05in
\proclaim{1.4. Proposition}Let $M$ be a separable {\rm II}$\sb 1$ factor with a Cartan subalgebra $A$. Then for 
any projection $p\in\Cal P(A)$ , the map
$\theta\in\text{Aut}_0(M;A)\rightarrow\theta_{\mid pMp}\in\text{Aut}_0(pMp;Ap)$ induces an isomorphism between the groups $\chi_0(M;A)$ and $\chi_0(pMp;Ap)$. In particular, $\chi_0(\sigma;G)$ is a stable orbit equivalence invariant.
 
 \endproclaim
\vskip .05in
{\it Proof.} Recall from [Po] that the restriction map $\theta\rightarrow\theta_{|pMp}$ induces an isomorphism between Out$_0(M;A)$ and Out$_0(pMp;Ap)$, which carries $\varepsilon_0(\overline{\text{Int}}_0(M;A))$ onto $\varepsilon_0(\overline{\text{Int}}_0(pMp;Ap))$. Lastly, recall that if $\theta$ is an automorphism of a II$_1$ factor $M$ which fixes a projection $p\in M$, then $\theta\in\text{Ct}(M)$ iff $\theta_{|pMp}\in\text{Ct}(pMp)$([Co1]).
\hfill $\blacksquare$
\vskip .2in

\head 2. Estimates of $\chi_0$.
\endhead
\vskip .1in
From this point on we will work with groups $G$ of the form $G=H\times K$, where $H$ is a non-amenable group and $K$ is an $\infty$ amenable group.

{\it 2.1. The model action.} Let $(X,\mu)$ be a standard probability space and let $\pi:H\rightarrow \text{Aut}(\Pi_{h\in H}(X,\mu)_h)$ be the left Bernoulli shift action given by $\pi_h((x_{h'})_{h'})=(x_{h^{-1}h'})_{h'},\forall h\in H,\forall x=(x_{h'})_{h'}\in\Pi_h (X,\mu)_h$.
Similarly, we let $\rho:K\rightarrow\text{Aut}(\prod_{k\in K}(X,\mu)_k)$ be the left Bernoulli shift action of $K$.

 Next, we define $\sigma:G\rightarrow\text{Aut}([\prod_{h\in H}(X,\mu)_h\times\prod_{k\in K}(X,\mu)_k])$ to be the product action given by $\sigma_{(h,k)}=\pi_h\times\rho_k,\forall h\in H, k\in K$. Note that $\sigma$ can also be viewed as the {\it generalized  Bernoulli action} ([Po3]) associated to the action of $G$ on $H\sqcup K$ given by $(h,k)\bullet x= hx,\forall x\in H$ and by $(h,k)\bullet x=kx,\forall x\in K$, $\forall g=(h,k)\in G=H\times K$.
Also, we denote $A=L^{\infty}(\prod_{h\in H}(X,\mu)_h)$, $B=L^{\infty}(\prod_{k\in K}(X,\mu)_k), C=A\overline{\otimes}B=L^{\infty}(\prod_{g\in H\sqcup K}(X,\mu)_g)$ and  $A_0=L^{\infty}((X,\mu)_e)\subset A$.
 
\vskip 0.05in
{\it 2.2. Quotient actions.} The actions of $G$ that we construct further arise as quotients of the {\bf model action $\sigma$}.
Let $\Gamma$ be a  group and take $\alpha:\Gamma\rightarrow\text{Aut}(A)$ and $\beta:\Gamma\rightarrow\text{Aut}(B)$  two free, m.p. actions of $\Gamma$ which commute with $
\pi$ and $\rho$, respectively.
Let $\delta$ denote the diagonal product action of $\Gamma$ on $C=A\overline{\otimes} 
B$ defined by $\delta(\gamma)=\alpha(\gamma)\otimes\beta(\gamma)$. 
If we set $C^\Gamma=\{x\in C\mid \delta(\gamma)(x)=x,  \forall \gamma\in\Gamma\}$ then, 
as $\delta$ and $\sigma$ commute, 
$G$  acts on $C^\Gamma$  and  this action, 
denoted $\sigma^\Gamma$, is   ergodic.
\vskip 0.05in

{\it Assumption.} Throughout this section we will assume that there is no $h\in H\setminus\{e\}$ such that $\pi(h)(a)=a,\forall a\in A^{\Gamma}$. Note that this is indeed  the case if $ A_0^{\Gamma}\not=\Bbb C1$.
\vskip  0.05in
Recall that the left Bernoulli shift action $\pi$  of a non-amenable group $H$ was shown to be strongly ergodic ([Sc2]), or equivalently, to satisfy $(A\rtimes_{\pi}H)'\cap A^{\omega}=\Bbb C1$, where $\omega$ is a free ultrafilter on $\Bbb N$([Po]). In what follows, we will need 
a more general statement:
\proclaim{2.3. Lemma} With the above notations we have that there exists $S\subset \Cal U(A^{\Gamma}\rtimes_{\pi^{\Gamma}}H)$ finite and $C>0$ such that $$\max_{u\in S}||[u,\xi]||_2\geq C||\xi-<{\xi},1>1||_2,\forall \xi\in L^2(A\rtimes_{\pi}H).$$In particular, $A\rtimes_{\pi}H$ is a  non $\Gamma$ type {\rm II}$\sb 1$ factor. 
\endproclaim
{\it Proof.} We begin by considering the representation of $H$ on $\Cal K=L^2(A\rtimes_{\pi} H)\ominus$ $l^2(H)$ given by $H\ni h\rightarrow\text{Ad}(u_h)\in\Cal U(\Cal K)$. If $\{\xi_0=1,\xi_1,..,\xi_n..\}$ is an orthonormal basis for $L^2(A_0)$, then $$\Cal B=\{ (\otimes_{h\in H}{\xi_{i_h}})u_{h'}|i_h\geq 0,\forall h\in H,1\leq |\{h|i_h\not=0\}|<\infty,h'\in H \}$$ gives an orthonormal basis for $\Cal K$. Since $H$ acts on $\Cal B$ and the stabilizers  $ H_{\xi}=\{ h\in H|\text{Ad}(u_h)(\xi)=\xi \}$ are finite  subgroups of $H$, $\forall\xi\in\Cal B$, we get that the representation of $H$ on $\Cal K$ is of the form 
$\bigoplus_i l^2(H/{H_i})$ for some finite subgroups $H_i$ of $H$.
Since $H$ is non-amenable, any representation of $H$ of the form  $\bigoplus_i l^2(H/{H_i})$, where $H_i$ are amenable subgroups of $H$, does not weakly contain the trivial representation ([Po]). Thus, there exist $k\in\Bbb N$, $h_1,h_2\dots h_k\in H$ and $c>0$ such that  $$\max_{i=\overline{1,k}} \| u_{h_i}xu_{h_i}^*-x\|_2\geq c\|x\|_2,\forall x\in \Cal K,$$ therefore, $$(*)\max_{i=\overline{1,k}} \| [u_{h_i},\xi]\|_2\geq c\|\xi-E_{L(H)}(\xi)\|_2,\forall \xi\in L^2(A\rtimes_{\pi} H).$$ 

Next, if $a\in \Cal U(A^{\Gamma})\setminus\Bbb C1$, then using the fact that Bernoulli shifts are mixing, we get that the set $F_a=\{h\in H|\pi_h(a)=a\}$ is finite and that $C_a=\inf_{h\in H\setminus F_a}||a-\pi_h(a)||_2>0$. If $x\in l^2(H)$, then $x=\Sigma_{h\in H}\tau(xu_{h^{-1}})u_h$ and we have that $$||[a,x]||_2^2=||\Sigma_{h\in H}\tau(xu_{h^{-1}})(a-\sigma_h(a))u_h||_2^2=$$ $$\Sigma_{h\in H}|\tau(xu_{h^{-1}})|^2 ||a-\sigma_h(a)||_2^2\geq C_a\Sigma_{h\in H\setminus F_a}|\tau(xu_{h^{-1}})|^2.$$
 Now, by our assumption, $\cap_{a\in\Cal U(A^{\Gamma})}F_a=\{h\in H|\sigma_h(a)=a,\forall a\in \Cal U(A^{\Gamma})\}=\{e\}$, thus   $\exists a_1,a_2,..,a_m\in\Cal U(A^{\Gamma})$ such that 
$\cap_{j=1}^{m}F_{a_i}=\{e\}$. By using the last inequality for $a\in\{a_1,a_2,..,a_m\}$ we get that $$(**)\max_{j=\overline{1,n}}||[a_j,x]||_2\geq  (\min_{j=\overline{1,m}}C_{a_j}/m)||x-\tau(x)||_2,\forall x\in l^2(H).$$

Combining $(*)$ and $(**)$, we get the conclusion for $S=\{u_{h_1},u_{h_2},..,u_{h_n},a_1,a_2,..,a_m\}$.
\hfill $\blacksquare$
\vskip 0.05in
{\it 2.4. Remark.} A consequence of the above proof  is that $L(H)'\cap L^2(A\rtimes_{\pi}H)\subset l^2(H)$.
We remark that for this  to be true we only need that  $\pi$ is weakly mixing (rather than the Bernoulli shift  of a non-amenable group). 
To see this, let $x=\Sigma_{h\in H}a_hu_h\in L(H)'\cap L^2(A\rtimes_{\pi}H)$, then $\pi_l(a_h)=a_{lhl^{-1}},\forall h,l\in H.$  In particular, this implies that $||a_h||_2=||a_{lhl^{-1}}||_2,\forall h,l\in H.$
For $h\in H\setminus \{e\}$ if the set $F=\{lhl^{-1}|l\in H\}$ is infinite then since $\Sigma_{h\in H}||a_h||_2^2=||x||_2^2<\infty$, we must have that $a_h=0$. On the other hand, if $F$ is finite then  $\Cal K=\text{span}\{a_h|h\in F\}\subset A$ is a finite dimensional vector space invaried under $\pi$. Since $\pi$ is weakly mixing, $K =\Bbb C1$, thus $a_h\in\Bbb C1$. 
\vskip 0.05in
The following result is well known, but we include a proof for the reader's convenience.
\proclaim{2.5. Lemma} Let $M$ and $N$ be two  {\rm II}$\sb 1$ factors and let $P$ be a subfactor of $M$.
 If there exists $S\subset\Cal U(P)$ finite and $C>0$ such that $$\max_{u\in S}||[u,\xi]||_2\geq C||\xi-<{\xi},1>1||_2,\forall \xi\in L^2M,$$ then 
 $P'\cap (M\overline{\otimes} N)^\omega= N^\omega$.
\endproclaim
\vskip 0.05in
{\it Proof.} Let $\{\eta_0=1,\eta_1,..,\eta_m...\}\subset L^2(N)$ be an orthonormal basis. Then any $x\in L^2(M\overline{\otimes} N)$ can be written as $x=\Sigma_{j\geq 0}\xi_j\otimes\eta_j$, for some $\xi_j\in L^2(M)$ with $\Sigma_{j\geq 0}||\xi_j||_2^2=||x||_2^2$.
Using the hypothesis we have that $$\Sigma_{u\in S}||[u,x]||_2^2=\Sigma_{u\in S}||\Sigma_{j\geq 0}[u,\xi_j]\otimes \eta_j||_2^2= $$ $$ \Sigma_{j\geq 0}\Sigma_{u\in S}||[u,\xi_j]||_2^2 \geq C^2\Sigma_{j\geq 0}||\xi_j-\tau(\xi_j)||_2^2=C^2||x-P_{L^2(N)}(x)||_2^2,$$ 

thus proving the conclusion.
\hfill $\blacksquare$
\vskip 0.07in

\proclaim{2.6. Corollary}With the notations we made before, we have that  $$(C^\Gamma\rtimes_{\sigma^\Gamma}G)'\cap( C^\Gamma\rtimes_{\sigma^\Gamma}G)^\omega\subset (B^\Gamma\rtimes_{\rho^\Gamma}K)^\omega.$$
\endproclaim
\vskip  .05in
{\it Proof.} By Lemma 2.3. we have that the inclusion $P=A^{\Gamma}\rtimes_{\pi^{\Gamma}}H\subset M=A\rtimes_{\pi}H$ satisfies the hypothesis of Lemma 2.5., thus  $$
(A^\Gamma\rtimes_{\pi^\Gamma}H)'\cap[(A\rtimes_{\pi}
H)\overline{\otimes}(B\rtimes_{\rho}K)]^{\omega}=(B\rtimes_{\rho}K)^{\omega}.$$ Since  $A^\Gamma\rtimes_{\pi^\Gamma}H\subset C^\Gamma\rtimes_{\sigma^\Gamma}G$ and $C^\Gamma\rtimes_{\sigma^\Gamma}G\subset [(A\rtimes_{\pi}H)\overline{\otimes}(B\rtimes_{\rho}K)]$, the corollary follows.
\hfill $\blacksquare$

\vskip 0.05in

From now on we  aim to calculate $\chi_0(\sigma^{\Gamma};G)$; to this end, we  give succesive estimates of $\chi_0$, which we improve as the class of groups $\Gamma$ and actions $\alpha,\beta$ that we consider, becomes more restrictive.
We start by showing  that the characters  of $\Gamma$ which have eigenvectors in both $\Cal U(A)$ and $\Cal U(B)$ naturally give rise to (non-trivial) elements of  $\chi_0(\sigma^{\Gamma};G)$. 
\vskip 0.05in
 
For an action $\delta$ of $\Gamma$ on $C$  and for a character $\eta\in$ Char($\Gamma$) we denote $\Cal {U}_\eta =\{a\in \Cal U(C) \mid
\delta(\gamma)(a)=\eta(\gamma)a , \forall \gamma\in\Gamma\}$ and we let $\text{Char}_{\delta}(\Gamma)=\{\eta\in \text{Char}(\Gamma)|\Cal U_{\eta}\not=\emptyset\}$([Po]).
Let $\eta\in\text{Char}_{\delta}(\Gamma)$, then for $a\in\Cal U_\eta$ the  1-cocycle   $w_g=a^*\sigma(a)$ is $C^\Gamma$-valued, thus
  $\theta=Ad(a)_{\mid C^\Gamma\rtimes_{\sigma^\Gamma}G}
\in$ Aut$\sb 0(C^\Gamma\rtimes_{\sigma^\Gamma}G;C^\Gamma)$.
Since the class of $\theta$ in Out$_0(C^\Gamma\rtimes_{\sigma^\Gamma}G;C^\Gamma)$
does not depend on the choice of $a\in\Cal U_\eta$, we get a well defined map 
$\phi:  
\text{Char}_{\delta}(\Gamma)\rightarrow$ H$^1(\sigma^\Gamma;G)$
  given by
$\phi(\eta)=\varepsilon_0(\theta)$.
\vskip  .1in
 
\proclaim{2.7. Proposition} For any $\eta\in \text{Char}_{\alpha}(\Gamma)\cap \text{Char}_{\beta}(\Gamma)$ we have that 
$\phi(\eta)\in\chi_0(\sigma^\Gamma;G)$. Moreover, the map $\phi:  
\text{Char}_{\alpha}(\Gamma)\cap \text{Char}_{\beta}(\Gamma)\rightarrow\chi_0(\sigma^\Gamma;G)$ is an injective group
homomorphism.
\endproclaim
 
{\it Proof.} For the first part  take $a_1\in \Cal U(A)\cap \Cal U_\eta$ and $a_2\in\Cal U(B)\cap U_\eta.$
Since $H$ is non-amenable, Corollary 2.3. implies that the central sequences of $C^\Gamma\rtimes_{\sigma^\Gamma}G$ asymptotically lie in $B^\Gamma\rtimes_{\rho^\Gamma}K$. Thus, as $a_1$ commutes with $B\rtimes_{\rho}K$, it results that 
$Ad(a_1)_{\mid C^\Gamma\rtimes_{\sigma^\Gamma}G}
\in$ Ct$(C^\Gamma\rtimes_{\sigma^\Gamma}G)$.  
Also, since $K$ is amenable we have that $\overline {\text {Int}_0}(B^{\Gamma}\rtimes_{\rho }K
;B^\Gamma)$=
Aut$_0(B^{\Gamma}\rtimes_{\rho} K;B^\Gamma)$([OW]) ,
  thus we deduce that
$Ad(a_2)_{\mid C^\Gamma\rtimes_{\sigma^\Gamma}G}
\in\overline {\text {Int}_0}(C^\Gamma\rtimes_{\sigma^\Gamma}G;C^\Gamma)$.
Combining the two conclusions we get the first assertion.
 
Finally, if we assume that $\phi$ is not injective then we get a non-trivial character of $\Gamma $, say $\eta$, and $a\in\Cal U_\eta$ which satisfies
  $Ad(a)_{\mid C^\Gamma\rtimes_{\sigma^\Gamma} G}
 \in$ Int$_0(C^\Gamma\rtimes_{\sigma^\Gamma}G;C^\Gamma)$. This
 implies that $\exists b\in\Cal U(C^\Gamma)$ such that $a^*\sigma(g)(a)=b^*\sigma(g)(b),
\forall g\in G$, but since $\sigma$ is ergodic we obtain that $a\in\Cal U(C^\Gamma)$, which leads to a contradiction, as $\eta$ is assumed non-trivial. 
\hfill $\blacksquare$
\vskip 0.1in

Next, we reduce the calculation of $\chi_0(\sigma;G)$ to the problem of finding the automorphisms of the hyperfinite II$_1$ factor $R$, which act trivially on certain central sequences of $R$.
\vskip 0.05in
{\it 2.8. Remark.}  To this aim, take 
$\theta\in\overline{\text {Int}_0}(\sigma^\Gamma;G)$, then $\theta=\lim_{n\rightarrow\infty}\text{Ad}(u_n)$ for some $u_n\in\Cal U(C^{\Gamma})$. 
Since by Corollary 2.3.,  $({C^\Gamma}\rtimes_{\sigma^\Gamma}G)'\cap {(C^{\Gamma})^{\omega}}\subset (B^{\Gamma})^{\omega}$, we deduce that there exist $a\in \Cal U(C^{\Gamma})$ and $a_n\in\Cal U(B^{\Gamma})$ such that $\theta=\text{Ad}(a)\lim_{n\rightarrow\infty}\text{Ad}(a_n)$([Co1],[Jo]). 
Thus, modulo $\text{Int}_0(\sigma^\Gamma;G)$, we can assume that $\theta=
\lim_{n\rightarrow\infty}\text{Ad}(a_n)$, for some $a_n\in\Cal U(B^{\Gamma})$, so in particular $\theta$ gives an automorphism of $B^\Gamma\rtimes_{\rho^\Gamma}K$.

Moreover, since $\theta$ is given by a
 $B^{\Gamma}$-valued 1-cocycle, hence $B$-valued, it extends to an automorphism, denoted $\tilde\theta$, of $R=B\rtimes_{\rho}K$. Since  $\tilde\theta\in{\text{Aut}}_0(\rho;K)$, we get a group morphism 
 $\psi:\overline{\text{Int}_0}(\sigma^{\Gamma};G)/{\text{Int}}_0(\sigma^{\Gamma};G)\rightarrow {\text{Aut}}_0(\rho;K)/{\text{Int}}_0(\rho;K)$ given by $\psi(\theta)=\tilde\theta$.
 
Now, remark that  since $R^{\Gamma}=B^{\Gamma}\rtimes_{\rho^{\Gamma}}K\subset C^{\Gamma}\rtimes_{\sigma^{\Gamma}}G$ and since $[R,A\rtimes_{\pi}H]=0$, we have the inclusion $R'\cap (R^{\Gamma})^{\omega}\subset (C^\Gamma\rtimes_{\sigma^\Gamma}G)'\cap( C^\Gamma\rtimes_{\sigma^\Gamma}G)^\omega$. Thus, if we further impose that $\theta\in\text{Ct}(C^{\Gamma}\rtimes_{\sigma^{\Gamma}}G)$, then $\tilde\theta$ acts trivially on $R'\cap (R^{\Gamma})^{\omega}$, i.e. $\tilde\theta\in\text{Ct}(R,R^{\Gamma})$.

We summarize the above discussion into the following:
\proclaim{2.9. Proposition}Assume that $\text{Char}_{\beta}(\Gamma)\subset \text{Char}_{\alpha}(\Gamma)$. Then the following is an exact sequence:
$$0\rightarrow\text{Char}_{\beta}(\Gamma)  
{\overset{\phi}\to\rightarrow}
\chi_0(\sigma^\Gamma;G){\overset{\psi}\to\rightarrow} 
({\text{Aut}}_0(\rho;K)\cap  \text{Ct}(R,R^{\Gamma})) /{\text{Int}}_0(\rho;K)$$
\endproclaim
{\it Proof.} To show that Ker$(\psi)\subset$ Ran$(\phi)$, let $\theta\in\chi_0(\sigma^{\Gamma};G)$ such that $\tilde\theta=\text{Ad}(a)$, for some $a\in\Cal U(B)$. But then, since $\theta(u_g)=a\sigma_g(a^*)\in C^{\Gamma},\forall g\in G$, we obtain that $\delta(\gamma)( a\sigma_g(a^*))=a\sigma_g(a^*),\forall g\in G,\forall\gamma\in\Gamma $. 
Using this relation and the fact that $\sigma$ and $\delta$ commute, we obtain that $$\sigma_g(a^*\delta(\gamma)(a))=a^*\delta(\gamma)(a),\forall g\in G,\forall\gamma\in\Gamma$$ and since $\sigma$ is ergodic, we conclude that $a^*\delta(\gamma(a))\in\Bbb C1,\forall \gamma\in\Gamma$. Thus, $\exists \eta\in\text{Char}_{\beta}(\Gamma)$ such that $\delta(\gamma)(a)=\eta(\gamma)a,  \forall \gamma\in\Gamma$, therefore $\theta\in\text{Ran}(\phi)$. The other inclusion follows easily, since  $\text{Ad}(a)\in \text{Ct}(R),\forall a\in\Cal U(B)$.
\hfill $\blacksquare$
\vskip .1in

We end this section by giving the motivation for the technical result that we prove in the next section. In view of Proposition 2.9.,  to compute $\chi_0(\sigma;G)$ it is natural to consider the following question: given a properly outer  action $\delta$ of a group $\Gamma$ on the hyperfinite II$_1$ factor $R$, describe   $\text{Ct}(R,R^{\Gamma})=\{\theta\in\text{Aut}(R)|\theta_{|R'\cap (R^{\Gamma})^{\omega}}=id\}$. 

If $\Gamma$ is finite, then  $\text{Ct}(R,R^{\Gamma})=\Gamma\text{Int}(R)$. Indeed, if we denote $R_{\omega}=R'\cap R^{\omega}$, then $R_{\omega}$ is II$_1$ factor on which $\Gamma$ and $\theta$ act properly ([Co1],[Jo]).  Since $\Gamma$ is finite we have that $(R^{\Gamma})^{\omega}=(R^{\omega})^{\Gamma}$(remark that this equality   is equivalent to the fact that the representation of $\Gamma$ on $L^2(R)\ominus L^2(R^{\Gamma})$  does not weakly contain the trivial representation, true if $\Gamma$ is a property (T) group), thus

$$R'\cap(R^{\Gamma})^{\omega}=R'\cap (R^{\omega})^{\Gamma}=(R'\cap R^{\omega})^{\Gamma}=R_{\omega}^{\Gamma},$$ hence   $\theta_{|R_{\omega}^{\Gamma}}=id _{|R_{\omega}^{\Gamma}}$.

 Recall now that if a finite group $\Gamma$ acts properly on a II$_1$ factor $M$ and $\theta\in\text{Aut}(M)$ is an automorphism which acts identically on $M^{\Gamma}$, then $\theta\in\Gamma$(note that
the last result also holds true if $\Gamma$ is a compact abelian group and $\text{Char}_{\delta}(\Gamma)=\text{Char}(\Gamma)$).
  Applying this to our situation, we can find $\gamma\in\Gamma$ such that $\theta_{|R_{\omega}}=\gamma _{|R_{\omega}}$, thus 
 ${\gamma}^{-1}\circ\theta\in\text{Ct}(R)=\text{Int}(R)$. As the other inclusion, $\Gamma\text{Int}(R)\subset\text{Ct}(R,R^{\Gamma})$, is always true, we get the desired equality for $\Gamma$ finite.  
Note that this equality already implies that $\chi_0(\sigma^{\Gamma};G)=\text{Char}(\Gamma)$ if $\Gamma$ is finite.

Turning to the case $\Gamma$ infinite, the above proof seems to fail as one would need the group $\Gamma$ to be both compact and rigid. However, in the next section  we show that  $\text{Ct}(R,R^{\Gamma})=\Gamma\text{Int}(R)$, given that the action $\Gamma\rightarrow\text{Aut}(R)$ is the "inductive limit" of actions of the form 
$\Gamma_n\rightarrow\text{Aut}(R_n)$, where $\Gamma_n$ are finite quotients of $\Gamma$ and  $R_n\subset R$ are subfactors exhausting $R$. 

\head 3.Main technical result. \endhead

\proclaim{3.1. Theorem} Let  $\Gamma=\bigoplus_{i\geq 0}\Delta_i$, where $(\Delta_i)_{i\geq 0}$ are finite groups and  let $\alpha:\Gamma\rightarrow \text{Aut}(R)$ be an action on the hyperfinite $II_1$ factor. Denote $\Gamma^n=\bigoplus_{i>n}\Delta_i$ and define $R^n=R^{\Gamma^n},\forall n\geq 0$.

(i) If  $\cup_{n\geq 0}R_n$ is weakly dense in $R$ , then $\alpha$ can be extended to a continuous action $\tilde\alpha:\tilde\Gamma=\prod_{i\geq 0}\Delta_i\rightarrow\text{Aut}(R)$.

(ii) Moreover, if  $({R^{\Gamma}})'\cap R=\Bbb C1$
 and if $\{\gamma\in\Gamma|\gamma_{|R_n}=\text{id}_{|R_n}\}=\Gamma^n,\forall n\geq 0$, then
$\text{Ct}(R,R^{\Gamma})=\tilde\Gamma\text{Int}(R)$.
 \endproclaim
{\it Proof.} (i) After identifying $\Gamma$ with $\alpha(\Gamma)\subset\text{Aut}(R)$, we need to prove that $\tilde\Gamma$ embedds into $\text{Aut}(R)$ in a continuous manner.
 Let $\gamma=(\delta_1,...,\delta_n,..)\in\tilde\Gamma$ and set $\gamma_n=(\delta_1,..,\delta_n,0,..,0,..)\in\Gamma$. We claim that $\exists\gamma:=\lim_{n\rightarrow\infty}\gamma_n$ in Aut($R$), or equivalently that $(\gamma_n(x))_n$ is a Cauchy sequence $\forall x\in R$. If $x\in\cup_{m\geq 0}R_m$, then the sequence  $(\gamma_n(x))_n$ eventually becomes constant, hence is Cauchy and as $\cup_{m\geq 0}R_m\subset R$ is a total set, the claim is proven. 

To prove the continuity of the embedding, let $\tilde{\Gamma}\ni\gamma_n=(\delta_i^n)_{i\geq 0}\rightarrow \gamma=(\delta_i)_{i\geq 0}\in\tilde{\Gamma}$ as $n\rightarrow\infty$.  This implies that $\forall i\geq 0,\exists N(i)\geq 0$ such that $\delta_i^n=\delta_i,\forall n\geq N(i)$. In turn, from this we deduce that if $x\in \cup_{m\geq 0}R_m$, then $\gamma_n(x)=\gamma(x),\forall n\geq N(x)$, thus $\gamma_n$ converges to $\gamma$ in the topology from $\text{Aut}(R)$.  
\vskip 0.05in (ii) 
For the second assertion, let $\theta\in\text{Ct}(R,R^{\Gamma})$, i.e. $\theta\in\text{Aut}(R)$ such that $\lim_{n\rightarrow\omega}||\theta(x_n)-x_n||_2=0,\forall (x_n)_n\in R'\cap(R^{\Gamma})^{\omega}$. Also, $\forall n\geq 0$, we denote $\Gamma_n=\bigoplus_{i\leq n}\Delta_i$.
\vskip 0.05in
{\it Step 1.}   $\exists N$ such that $\forall n\geq N$ we have that $||\theta(x)-x||_2\leq1/2,\forall x\in \Cal U(P_n)$, where $P_n:=R_n'\cap (R_n^{\Gamma_n})^{\omega}={(R_n'\cap R_n^{\omega})}^{\Gamma_n}$.
\vskip 0.05in

 This follows  from the following:
\proclaim{3.2. Lemma} Let $M$ be a II$_1$ factor and let $N,\{M_n\}_{n\geq 0}\subset M$ be  von Neumann subalgebras  such that $\lim_{n\rightarrow\infty}||x-E_{M_n}(x)||_2=0,\forall x\in M$.  If $\theta\in\text{Aut}(M)$  acts trivially on $M'\cap N^{\omega}$, then $\forall \varepsilon>0,\exists N=N(\varepsilon)\in\Bbb N$ such that $||\theta(u)-u||_2\leq \varepsilon,\forall u\in\Cal U(M_n' \cap N^{\omega}),\forall n\geq N(\varepsilon)$. 
\endproclaim
{\it Proof of lemma 3.2.} If the conclusion fails, then $\exists\varepsilon>0$, an increasing subsequence $\{k_n\}$ of $\Bbb N$ and $u_n\in\Cal U(M_{k_n}'\cap N^{\omega})$ such that $||\theta(u_n)-u_n||_2> \varepsilon,\forall n\in\Bbb N$.
Since $\theta$ acts trivially on $M'\cap N^{\omega}$, there exists $F\subset M$ finite and $\delta>0$ such that if $x\in (N)_1$ satisfies $||[x,y]||_2\leq\delta,\forall y\in F$, then $||\theta(x)-x||_2\leq \varepsilon/2.$

Representing $u_n$ as $(u_n^m)_m$, where $u_n^m\in\Cal U(N)$ and using the first inequality, we deduce that $\exists m_n\in\Bbb N$ such that $v_n=u_n^{m_n}$ satisfies $||[v_n,E_{M_{k_n}}(x)]||_2\leq\delta/2,\forall x\in F$ and $||\theta(v_n)-v_n||>\varepsilon,\forall n\in\Bbb N$. Finally, since $\lim_{n\rightarrow\infty}||E_{M_{k_n}}(x)-x||_2=0,\forall x\in F$ and since $||[v_n,x]]||_2\leq ||[v_n,E_{M_{k_n}}(x)]||_2+2||E_{M_{k_n}}(x)-x||_2\leq \delta/2+2||E_{M_{k_n}}(x)-x||_2$, we get a contradiction.
\hfill $\blacksquare$
\vskip 0.05in
Going back to the proof of Theorem 3.1., we can apply  Lemma 3.2. to find $N=N(1/2)$ such that  $$||\theta(x)-x||_2\leq1/2,\forall x\in \Cal U(R_n'\cap (R^{\Gamma})^{\omega}),\forall n\geq N.$$ Since by  the definition of $R_n$ we have that $R_n^{\Gamma}=R_n^{\Gamma_n}
$, we deduce that $P_n\subset R_n'\cap (R^{\Gamma})^{\omega}$, thus proving this step.

Next, by standard averaging techniques, we get :
\vskip 0.05in
{\it Step 2.} $\forall n\geq N, \exists w_n\in R^{\omega},w_n\not=0$ such that $\theta(x)w_n=w_nx,\forall x\in \Cal U(P_n)$.
\vskip 0.05in
Indeed, from the inequality $||\theta(x)-x||_2\leq 1/2,\forall x\in\Cal U(P_n)$, we deduce that the minimal $||.||_2$ element of the convex set $\Cal K_n=\overline{co}^{w}\{\theta(x)x^*|x\in\Cal U(P_n)\}\subset (R^{\omega})_1$ satisfies $||w_n-1||_2\leq 1/2$, thus $w_n\not=0$. Moreover, by the uniqueness of $w_n$ and by the fact that $\Cal K_n$ is invaried under the transformations $u\rightarrow \theta(x)ux^*,\forall x\in\Cal U(P_n)$, it follows that $w_n$  satisfies the required identity.
\vskip 0.05in
{\it Step 3.} 
 \vskip 0.05in

\proclaim{3.3. Lemma}  Let $\Gamma$ be a finite group acting on the inclusion $N\subset M$ of two finite von Neumann algebras such that $N'\cap (M\rtimes\Gamma)\subset M$. 
If  $\theta\in \text{Aut}(M)$ satisfies
$\theta(x)w=wx,\forall x\in N^{\Gamma}$, for some non-zero element $w\in M$, then $\exists v\in M,v\not=0$ and $\gamma\in\Gamma$ such that $\theta(x)v=v\gamma(x),\forall x\in N$.
\endproclaim
{\it Proof of lemma 3.3.} 
If $x\in\Cal  U(N)$, then $\Sigma_{\gamma\in\Gamma}\gamma(x)\in N^{\Gamma}$, thus $$\Sigma_{\gamma\in\Gamma}\theta(\gamma(x))w=w\Sigma_{\gamma\in\Gamma}\gamma(x), \forall x\in \Cal U(N).$$ Using this relation  we get that
 $$\Sigma_{\gamma\in\Gamma}\theta(\gamma(x))wx^*=w\Sigma_{\gamma\in\Gamma}\gamma(x)x^*,
\forall x\in \Cal U(N).$$
Let $\tilde M:=\bigoplus_{\gamma\in\Gamma}(M)_{\gamma}$ and define the following two uniformly bounded convex subsets of $\tilde M$:
$$\Cal K_1=\overline{co}^{w}\{\oplus_{\gamma\in\Gamma}(\theta(\gamma(x))wx^*)|x\in\Cal U(N)\},$$
$$\Cal K_2=\overline{co}^{w}\{\oplus_{\gamma\in\Gamma}(\gamma(x)x^*)|x\in\Cal U(N)\}.$$
Then the above relation translates into:  $$\Sigma_{\gamma}x_{\gamma}^1=\Sigma_{\gamma}wx_{\gamma}^2,\forall x^i=\oplus _{\gamma\in\Gamma}x_{\gamma}^i\in\Cal K_i, i=1,2.$$
Let $x^i=\oplus _{\gamma\in\Gamma}x_{\gamma}^i\in\Cal K_i, i=1,2$ be the elements of minimal $||.||_2$. Then note that $\forall u\in\Cal U(N)$, $\Cal K_1$ and $\Cal K_2$ are invaried under the transformations $\oplus_{\gamma}x_{\gamma}^1\rightarrow \oplus_{\gamma}\theta
(\gamma(u))x_{\gamma}^1u^*$ and $\oplus_{\gamma}x_{\gamma}^2\rightarrow\oplus_{\gamma}\gamma(u)x_{\gamma}^2u^*$, respectively.

Thus, since these transformations are norm preserving, it follows by the uniqueness of $x_1$ and $x_2$ that $\theta
(\gamma(u))x_{\gamma}^1u^*=x_{\gamma}^1$ and $\gamma(u)x_{\gamma}^2u^*=x_{\gamma}^2,\forall \gamma\in\Gamma,\forall u\in\Cal U(N)$. Now, note that the condition $N'\cap (M\rtimes\Gamma)\subset M$ is equivalent to the following:
  if $ \gamma\in\Gamma\setminus\{e\}$ and  $w\in M$, then   $$w\gamma(x)=xw,\forall x\in N\Rightarrow w=0.$$
Thus $x_{\gamma}^2=0,\forall\gamma\in\Gamma\setminus\{e\}$ and since $x_{e}^2=1$ and $\Sigma_{\gamma}x_{\gamma}^1=\Sigma_{\gamma}wx_{\gamma}^2$, it follows that there exists $\gamma\in\Gamma$ such that $x_{\gamma}^1\not=0$.
\hfill $\blacksquare$

\vskip 0.05in

{\it Step 4.} $\forall n\geq N,\exists v_n\in R^{\omega},v_n\not= 0$ and $\gamma_n\in \Gamma_n$ such that $\theta(x)v_n=v_n\gamma_n(x),\forall x\in S_n:=R_n'\cap {R_n^{\omega}}$. 
\vskip 0.05in
Note first that $P_n=R_n'\cap (R_n^{\Gamma_n})^{\omega}=S_n^{\Gamma_n}$ and that $\Gamma_n$ acts on the inclusion $S_n\subset R^{\omega}$. Recall that by {\it Step 2} we have that $\theta(x)w_n=w_nx,\forall x\in S_n^{\Gamma_n}$.  
Thus, provided that we can show that there is no non-zero $w\in R^{\omega}$ such that $\gamma(x)w=wx,\forall x\in S_n=R_n'\cap R_n^{\omega}$ for some $\gamma\in\Gamma_n\setminus\{e\}$, the conclusion of Lemma 3.3. gives our claim. 

To disprove the last equality, we show the following lemma, which we will also use subsequently.
\vskip 0.05in
\proclaim{3.4. Lemma}Let $M$ be a $II_1$ factor and let $R$ be an irreducible hyperfinite subfactor. If $\alpha\in\text{Aut}(M)$ is an automorphism such that $\exists w\in M^{\omega},w\not=0, \alpha(x)w=wx,\forall x\in R'\cap R^{\omega}$, then there exists $v\in\Cal U(M)$  such that $\alpha(x)=vxv^*,\forall x\in R.$
\endproclaim
{\it Proof of lemma 3.4.} Since $R$ is hyperfinite, we can decompose $R=\overline{\otimes}_{k\in\Bbb N}(\Bbb M_{2\times 2}(\Bbb C))_k$.  Denote $R_{m}=\overline{\otimes}_{k\geq m}(\Bbb M_{2\times 2}(\Bbb C))_k$ and represent $w=(w_n)_n\in M^{\omega}$. 
We claim that $\exists m\in \Bbb N$ such that $$||\alpha(x)w_m-w_mx||_2\leq ||w_m||_2,\forall x\in\Cal U(R_m).$$ Indeed, because if for every $m\in\Bbb N$ we can find $x_m\in \Cal U(R_m)$ such that $||\alpha(x_m)w_m-w_mx_m||_2> ||w_m||_2$, then $x=(x_m)_m\in \Cal U(R'\cap R^{\omega})$ and $ ||\alpha(x)w-wx||_2=\lim_{n\rightarrow\omega}||\alpha(x_m)w_m-w_mx_m||_2\geq\lim_{n\rightarrow \omega}||w_m||_2=||w||_2$, a contradiction.

Next, a simple computation shows that $$\text{Re} {\tau(w_m^*\alpha(x)w_mx^*)}\geq ||w_m||_2/2,\forall x\in \Cal U(R_m).$$
Thus, if $v$ denotes the element of minimal $||.||_2$  in $C=\overline{co}^{w}\{\alpha(x)w_mx^*|x\in\Cal U(R_m)\}$, then Re$\tau(w_m^*v)\geq||w_m||_2/2$, hence $v\not=0$ and $\alpha(x)v=vx,\forall x\in\Cal U(R_m)$. 

Since $R'\cap M=\Bbb C1$, we get that $R_m'\cap M=\overline{\otimes}_{k<m}(\Bbb M_{2\times 2}(\Bbb C))_k$ (this is true since if $A\subset B$ are $C^*$-algebras, then ${A^{(n)}}'\cap \Bbb M_{ n}(B)=\Bbb M_{ n}(A'\cap B)$). Also, because $v^*v\in R_m'\cap M,vv^*\in \alpha(R_m)'\cap M$, we deduce that there exists $p,q\in\Cal P(R_m'\cap M)\subset \Cal P(R)$ and $u\in\Cal U(M)$ such that $\alpha(x\otimes p)=u(x\otimes q)u^*,\forall x\in R_m$. Finally, since $p,q$ have central supports equal to 1 in  $\overline{\otimes}_{k<m}(\Bbb M_{2\times 2}(\Bbb C))_k$, we conclude that $\exists v\in\Cal U(M)$ such that $\alpha(x)v=vx,\forall x\in R$.
\hfill $\blacksquare$
\vskip 0.05in
Assume by contradiction that $\exists w\in R^{\omega}$ such that $\gamma(x)w=wx,\forall x\in S_n=R_n'\cap R_n^{\omega}$ for some $\gamma\in\Gamma_n\setminus\{e\}$. Since $R^{\Gamma}\subset R_n$ and $(R^{\Gamma})'\cap R=\Bbb C1$, we get that $R_n'\cap R=\Bbb C1$, thus Lemma 3.4. implies that $\exists v\in\Cal U(R),$ $\gamma(x)=vxv^*,\forall x\in R_n.$

Using again the inclusion $R^{\Gamma}\subset R_n$ , we infer that $v\in {R^{\Gamma}}'\cap R=\Bbb C1$, thus, $\gamma(x)=x,\forall x\in R_n$.
Now, the hypothesis  entails that $\gamma\in \Gamma^n$, therefore $\gamma\in\Gamma^n\cap\Gamma_n=\{e\}$, contradiction. 
\vskip 0.05in
{\it Step 5.} $\forall n\geq N,\exists u_n\in \Cal U(R)$ and $\gamma_n\in \Gamma_n$ such that $\theta(x)u_n=u_n\gamma_n(x),\forall x\in R_n$. 
\vskip 0.05in
This step follows directly from {\it Step 4} and Lemma 3.4. by using the fact that $R_n$ has trivial relative commutant in $R,\forall n\geq 0$.
\vskip 0.05in
{\it Step 6.} $\theta\in\tilde{\Gamma}\text{Int}(R)$.
\vskip 0.05in
Following {\it Step 5}, $\forall n\geq N$, we can find $u_n\in\Cal U(R),\gamma_n\in\Gamma_n$ such that $\theta(x)=u_n\gamma_n(x)u_n^*,\forall x\in R_n.$
Write $\gamma_{n+1}=\gamma_{n}'\delta_{n+1}$, where $\gamma_{n}'\in\Gamma_n,\delta_{n+1}\in\Delta_{n+1},\forall n\geq N $ and let $\gamma_{N}=(\delta_0,\delta_1,..,\delta_N,0,..,0..)\in\Gamma_N$. Using the fact that $$u_n\gamma_n(x)u_n^*=\theta(x)=u_{n+1}\gamma_{n+1}(x)u_{n+1}^*,\forall x\in R_n$$ and that $\delta_{n+1}$ acts identically on $R_n$, we get that $$\gamma_n\circ{\gamma_n'}^{-1}(x)=(u_n^*u_{n+1})x({u_{n+1}^*}u_n),\forall x\in R_n.$$

Since $\gamma_n\circ{\gamma_n'}^{-1}\in\Gamma_n$, by reasoning as in the proof of {\it Step 4} we deduce that $\gamma_n=\gamma_{n+1}'$ and that $u_n^*u_{n+1}\in\Bbb C1$.
 Consequently, if we let $\gamma=(\delta_n)_{n\geq 0}\in\tilde\Gamma$ and $u=u_N$, then $\theta(x)=u\gamma(x)u^*,\forall x\in\cup_{n\geq N}R_n$, thus finishing the proof. 
\hfill $\blacksquare$

\head 4.Calculation of $\chi_0(\sigma^{\Gamma};G)$. \endhead
In this section we use the results of the previous two sections to get concrete computations of $\chi_0$. 
As before, let $H$ be a non-amenable group and $K$ be an $\infty$ amenable group. Also, we fix $(X_0,\mu_0)$, a standard probability space and we define $$(X,\mu)=\prod_{i\geq 0}(X_0,\mu_0)_i.$$
Then,  as in  Section 2, we let $A=L^{\infty}(\prod_{h\in H}(X,\mu)_h)$, $B=L^{\infty}(\prod_{k\in K}(X,\mu)_k), C=A\overline{\otimes}B=L^{\infty}(\prod_{g\in H\sqcup K}(X,\mu)_g)$ and we denote  by $\pi,\rho$ the left Bernoulli shift actions of $H,K$ on $A,B$, respectively and by $\sigma$ the product action of $G=H\times K$ on $C$.

Let $\{\Delta_i\}_{i\geq 0}$ be non-trivial finite groups and denote $\Gamma=\bigoplus_{i\geq 0}\Delta_i$, $\Gamma^n=\bigoplus_{i>n}\Delta_i,\Gamma_n=\bigoplus_{i\leq n}\Delta_i, \forall n\geq 0.$ 
For $i\geq 0$, let $\alpha_i:\Delta_i\rightarrow\text{Aut}(X_0,\mu_0)$ be a free action. 
Then, we define  $\alpha:\Gamma\rightarrow\text{Aut}(\prod_{h\in H}(X,\mu)_h=\prod_{i\geq 0,h\in H}(X_0,\mu_0)_{i,h})$ to be  given by $$\alpha(\gamma)((x_{i,h})_{i,h})=((\alpha_i(\delta_i)(x_{i,h}))_{i,h},$$
 $\forall x={(x_{i,h})}_{i,h}\in \prod_{i\geq 0,h\in H}(X_0,\mu_0)_{i,h},\forall\gamma=(\delta_i)_{i\geq 0}\in\bigoplus_{i\geq 0}\Delta_i=\Gamma.$ Similarly, we let $\beta:\Gamma\rightarrow\text{Aut}(\prod_{k\in K}(X,\mu)_k=\prod_{i\geq 0,k\in K}(X_0,\mu_0)_{i,k})$. We  also denote by $\alpha,\beta$ the induced actions on $A,B$.

Now, consider the diagonal product action $\delta:\Gamma\rightarrow\text{Aut}(C)$ given by $$\delta(\gamma)=\alpha(\gamma)\otimes\beta(\gamma),\forall\gamma\in\Gamma.$$
Since $[\alpha,\pi]=0$ and $[\beta,\rho]=0$, it follows that $[\delta,\sigma]=0$, thus $G$ acts on $C^{\Gamma}=\{x\in C|\delta(\gamma)(x)=x,\forall \gamma\in\Gamma\}$, hence we can define $\sigma^{\Gamma}=\sigma_{|C^{\Gamma}}.$

\proclaim {4.1. Theorem}  $\sigma^{\Gamma}:G\rightarrow\text{Aut}(C^{\Gamma})$ is a free, ergodic, integral preserving action and $$\phi:\text{Char}(\Gamma)\simeq \chi_0(\sigma^{\Gamma};G)$$
\endproclaim
{\it Proof.} If we let $R=B\rtimes_{\rho}K$, then since $K$ is amenable and $\rho$ is ergodic, $R$ is the hyperfinite II$_1$ factor ([OW]).  Moreover, since $[\beta,\rho]=0$, we have that $\Gamma$ acts on $R$ by $\beta(\gamma)(\Sigma_{k\in K}a_ku_k)=\Sigma_{k\in K}\beta(a_k)u_k$ and that  $R^{\Gamma}=B^{\Gamma}\rtimes_{\rho^{\Gamma}}K$.
\vskip 0.05in

{\it Claim 1.} If $\theta\in\text{Aut}(R)$ acts identically on $R'\cap (R^{\Gamma})^{\omega}$, then $\theta\in\tilde\Gamma\text{Int}(R)$.
\vskip 0.05in

To prove this it is sufficient to verify  the conditions  of Theorem 3.1.
Note that $\Gamma^n$ acts trivially on $\prod_{n\geq i\geq 0,k\in K}(X_0,\mu_0)_{i,k}$ and also on  $K$, thus implying that $$L^{\infty}(\prod_{n\geq i\geq 0,k\in K}(X_0,\mu_0)_{(i,k)})\rtimes K\subset R_n=R^{\Gamma^n},\forall n\geq 0.$$ Further, this clearly implies that $\overline{\cup_{n\geq 0}R_n}^{w}=R$.
Next, Lemma 2.1. implies that $(B^{\Gamma}\rtimes_{\rho^{\Gamma}}K)'\cap (B\rtimes_{\rho}K) =\Bbb C1$  or, equivalently, that ${R^{\Gamma}}'\cap R=\Bbb C1$. 
Finally, if $\gamma=(\delta_1,..,\delta_n,0,..)\in\Gamma_n,$  where $\delta_i\in\Delta_i$, acts trivially on $R_n$, then in particular, it acts trivially on $\prod_{n\geq i\geq 0,k\in K}(X_0,\mu_0)_{i,k}$.
 This in turn forces $\alpha_i(\delta_i)=$id, hence $\delta_i=e,\forall i=\overline{1,n}$, as $\alpha_i$ are assumed free.

\vskip 0.05in
{\it Claim 2.} If  $g=(h,k)\in G\setminus\{e\}$  then there exists no non-zero $p\in C$ such that $\sigma_g(c)p=cp,\forall c\in A^{\Gamma}\overline{\otimes}B^{\Gamma}\subset C^{\Gamma}$. 
\vskip 0.05in

To see this, let $a\in L^{\infty}(\prod_{i\geq 0}(X_0,\mu_0)_{i,e})^{\Gamma}\subset A^{\Gamma}$ and $b\in L^{\infty}(\prod_{i\geq 0}(X_0,\mu_0)_{i,e})^{\Gamma}\subset B^{\Gamma}$, such that $\tau(a)=\tau(b)=0$. Thus $\tau(a^*\pi_h(a))=\tau(b^*\rho_k(b))=0,\forall h\in H\setminus\{e\},k\in K\setminus\{k\}.$ 

Recall from [Po1] that Bernoulli shift actions are 2-mixing in the following sense:

$$\lim_{h_1,h_2\rightarrow\infty}|\tau(x_0\pi_{h_1}(x_1)\pi_{h_2}(x_2))-\tau(x_0)\tau(\pi_{h_1}(x_1)\pi_{h_2}(x_2))|=0,\forall x_0,x_1,x_2\in A.$$

 Now, let $h_n\in H,k_n\in K$, both going to $\infty$, and define $c_n=\pi_{h_n}(a)\otimes\rho_{k_n}(b)\in A^{\Gamma}\otimes B^{\Gamma}$. Fix $x\in A,y\in B$ and denote $z=x\otimes y$, then 
$$\tau(c_n^*\sigma_g(c_n)z)=\tau(\pi_{h_n}(a)^*\pi_{h}(\pi_{h_n}(a))x)\tau(\rho_{k_n}(b)^*\rho_{k}(\rho_{k_n}(b))y).$$

Using the 2-mixingness, we deduce that $$\lim_{n\rightarrow\infty} |\tau(c_n^*\sigma_g(c_n)z)-\tau(x)\tau(y)\tau(\pi_{h_n}(a)^*\pi_{hh_n}(a))\tau(\rho_{k_n}(b)^*\rho_{kk_n}(b))|=0.$$

On the other hand, $$\tau(\pi_{h_n}(a)^*\pi_{hh_n}(a))\tau(\rho_{k_n}(b)^*\rho_{kk_n}(b))=\delta_{h_n,hh_n}\delta_{k_n,kk_n}||a||_2^2||b||_2^2=\delta_{(h,k),(e,e)}||a||_2^2||b||_2^2=0,$$
thus $\lim\tau(c_n^*\sigma_g(c_n)z)=0,\forall z\in A\otimes B$ and, as $A\otimes B$ is $||.||_2$-dense in $C=A\overline{\otimes}B$, the claim follows.

\vskip 0.06in
 
Note that the above claim implies directly that $\sigma^{\Gamma}$ is a free action, while ergodicity follows from that of $\sigma$, which in turn is implied by the fact that both $\pi$ and $\rho$ are weakly mixing.
 
\vskip 0.05in
 {\it Claim 3.} $\tilde\Gamma\text{Int}(R)\cap \text{Aut}_0(R;B)=\text{Int}_0(R;B).$
\vskip 0.05in

Indeed, if $\tilde\gamma\in\tilde\Gamma\setminus\{e\}$ and $u\in\Cal U(R)$ are such that $\theta=\text{Ad}(u)\circ\tilde\gamma\in\text{Aut}_0(R;B)$, then $\tilde\gamma(x)=\text{Ad}(u^*)(x),\forall x\in B$. In particular, this implies that $u\in\Cal N_{R}(B)$, thus we can find a non-zero projection $p\in B$ and $k\in K\setminus\{e\}$ such that $\rho_k(b)p=Ad(u^*)(b)p=\tilde\gamma(b)p,\forall x\in B$, which entails that $\rho_k(b)p=bp,\forall b\in B^{\Gamma}.$ 

Denote $g=(e,k)\in H\times K=G$, then  we have that  $\sigma_g(a\otimes b)(1\otimes p)=(a\otimes b)(1\otimes p),\forall a\in A^{\Gamma},\forall b\in B^{\Gamma}$, thus by the second claim, $g=e$, hence $k=e$, a contradiction.  
\vskip 0.05in

Combining {\it Claim 1}  and {\it Claim 3}, we get that $$\text{Ct}(R,R^{\Gamma})\cap \text{Aut}_0(R;B)={\tilde\Gamma}\text{Int}(R)\cap{\text{Aut}}_0(R;B)=\text{Int}_0(R;B),$$

 thus, by Proposition 2.9.,  $\phi:{\text{Char}}_{\beta}(\Gamma)\rightarrow\chi_0(\sigma^{\Gamma};G)$ is an isomorphism. Finally, it is clear that ${\text{Char}}_{\beta}(\Gamma)=\text{Char}(\Gamma)$, which finishes the proof.
\hfill $\blacksquare$
\proclaim {4.2. Corollary} Let $G=H\times K$, where $H$ is a non-amenable  group and $K$ is an $\infty$ amenable group. Then, for any group $\Lambda=\Pi_{i\geq 0}\Lambda_i$, where $\Lambda_i$ are finite abelian groups, we can find a free, ergodic, measure preserving action $\sigma_{\Lambda}$ of $G$ on a standard probability space such that $\chi_0(\sigma_{\Lambda};G)=\Lambda$. Thus, any such group $G$ admits uncountably many non stably orbit equivalent actions. 
\endproclaim 
\vskip 0.05in
{\it 4.3. Final remarks.}
 (i). Similarly, it can  be shown that a non-amenable group $G$ which is the infinite sum of non-trivial groups has   uncountably many non-OE actions (e.g. $G=\bigoplus_{i\geq 0}\Bbb F_n,\forall n\geq 2$). 

(ii). A related OE  invariant for actions $\sigma:G\rightarrow$ Aut($X,\mu$) that one can  consider is
given by $\overline{\text B^1}(\sigma,G)$/ B$^1(\sigma,G)\simeq\overline{{\text {Int}}_0}(M;A)/$ Int$_0(M;A)$. 
Note that this invariant is a subgroup of the 1-cohomology group H$^1(\sigma;G)$ which is isomorphic
to Kawahigashi's relative $\chi$ invariant of the corresponding Cartan subalgebra inclusion, $\chi(M,A)$ ([Ka]).
While it seems difficult to calculate such an invariant (as it is an infinite, non-Polish group) we note
that  if a group $\Cal G$ is isomorphic to $\overline{\text B^1}(\sigma,G)$/ B$^1(\sigma,G)$, for some countable, discrete group $G$, then we can replace $G$ by any non-amenable group which admits $G$ as a quotient (e.g. $\Bbb F_{\infty}$).

(iii). Note that the results that we prove also work if we consider actions on the hyperfinite II$_1$ factor $R$ (instead of a diffuse abelian von Neumann algebra), thus rendering the analogous conclusion: any group $G$ as in the above Corollary has uncountably many non-cocycle conjugated actions on $R$. The only difference is the following: in the proof of Proposition 2.7., we need to use a result of Ocneanu ([Oc]) stating that 1-cocycles for actions of amenable groups on the hyperfinite II$_1$ factor are approximately 1-coboundaries, instead of the analogous result for actions on a probability space([OW]).
 
\vskip .1in
\head  References\endhead
 
\item {[Co1]} A. Connes: {\it Sur la classification des facteurs de type \rm II}, C. R. Acad. Sci. Paris SŽr.I Math.
{\bf 281} (1975), A13ÐA15. 
\item {[Co2]} A. Connes: {\it Classification of injective factors}, Ann. of Math., {\bf 104}
(1976), 73-115.
\item {[CFW]} A. Connes, J. Feldman, B. Weiss: {\it An amenable equivalence relation is generated 
a single transformation}, Ergodic Theory Dynamical systems {\bf 1} (1981), 431-450.
\item {[CW]} A.Connes, B. Weiss: {\it Property (T) and asymptotically invariant sequences}, 
Israel J. Math. {\bf 37} (1980),  209-210. 
\item {[Dy]} H.Dye: {\it On groups of measure preserving transformations, \rm I}, Amer. J. Math.
{\bf 81} (1959), 119-159. 
\item {[FM]} J. Feldman, C.C. Moore: {\it Ergodic equivalence relations,  cohomology, and von Neumann
 algebras \rm II}, Trans. Amer. Math. Soc. {\bf 234} (1977), 325-359.
\item {[GaPo]} D. Gaboriau, S. Popa: {\it An Uncountable Family of Non Orbit Equivalent Actions of $\Bbb F_n$}
J. Amer. Math. Soc. 18 (2005), no. 3, 547-559. 
 \item {[Hj]} G. Hjorth: {\it A converse to Dye's Theorem}, Trans. Amer. Math. Soc. {\bf 357} (2005), 3083-3103.
 \item {[Ka]} Y. Kawahigashi: {\it Centrally trivial automorphisms and an analogue of Connes' $\chi(M)$
for subfactors}, Duke Math. J. {\bf 71} (1993), 93-118.
\item {[Jo]} V.F.R. Jones: {\it Notes on Connes' invariant $\chi$(M)}, unpublished.
\item {[MoSh]} N. Monod, Y. Shalom: {\it Orbit equivalence rigidity and bounded cohomology}, Preprint 2002, to appear in Ann. of Math.
\item {[MvN]} F. Murray, J. von Neumann: {\it Rings of operators IV}, Ann. Math. {\bf 44} (1943), 716-808.	
\item {[Oc]} A. Ocneanu: {\it Actions of discrete amenable groups on von Neumann algebras}, Lecture Notes in Mathematics, 1138, Springer Verlag, Berlin, 1985.
\item {[OW]} D.Ornstein, B.Weiss: {\it Ergodic theory of amenable groups. I. The Rokhlin lemma.}, Bull. Amer. Math. Soc. (N.S.) {\bf 1} (1980), 161-164.
\item {[Po]} S. Popa: {\it Some computations of 1-cohomology groups and construction of non orbit
equivalent actions}, math. OA/0407199, to appear in Journal de l'Inst. de Math. de Jussieu.
\item {[Po1]} S. Popa: {\it Strong rigidity of II$_1$ factors arising from malleable actions of $w$-rigid groups, II}, math. OA/0407103.

\item {[Po2]} S.Popa: {\it Some rigidity results for non-commutative Bernoulli shifts}, Journal of Functional Analysis {\bf 230} (2006), 273-328
\item {[PoSa]} S. Popa,R. Sasyk: {\it On the cohomology of actions of groups by Bernoulii shifts},
math. OA/ 0310211.   
\item {[Si]} I.M. Singer: {\it Automorphisms of finite factors}, Amer. J. Math. {\bf 77} (1955), 117-133.
\item {[Sc1]} K. Schmidt: {\it Asymptotically invariant sequences and an action of SL(2,$\Bbb Z$) on the
2-sphere}, Israel J. Math. {\bf 37} (1980), 193-208.
\item {[Sc2]} K.  Schmidt: {\it Amenability, Kazhdan's property T, strong ergodicity and invariant
means for ergodic group-actions,} Ergodic. Th.  Dynam. Sys. {\bf 1} (1981),
 223-236.

\enddocument